\documentstyle[11pt]{article}

\newtheorem{theo}{Theorem}
\newtheorem{coro}{Corollary}
\newtheorem{prop}{Proposition}
\newtheorem{lemm}{Lemma}

\begin{document}

\def\ot{\otimes}
\def\we{\wedge}
\def\wec{\wedge\cdots\wedge}
\def\op{\oplus}
\def\ra{\rightarrow}
\def\lra{\longrightarrow}
\def\fso{\mathfrak so}
\def\cO{\mathcal{O}}
\def\cS{\mathcal{S}}
\def\fsl{\mathfrak{sl}}
\def\PP{\mathbb P}\def\PP{\mathbb P}\def\ZZ{\mathbb Z}\def\CC{\mathbb C}
\def\RR{\mathbb R}\def\HH{\mathbb H}\def\OO{\mathbb O}

\title{On the asymptotics of Kronecker coefficients, 2}
\author{Laurent \sc{Manivel} \\  \\
UMI 3457 CNRS/Centre de Recherches Math\'ematiques,\\
Universit\'e de Montr\'eal, Canada\\
\texttt{manivel@math.cnrs.fr} }
\maketitle

\begin{abstract}
Kronecker coefficients encode the tensor products of complex irreducible representations of symmetric groups. Their stability properties have been considered recently by several authors (Vallejo, Pak and Panova, Stembridge). In \cite{ak1} we 
described a geometric method, based on Schur-Weyl duality, that allows to produce huge series of instances of this phenomenon.
In this note we show how to go beyond these so-called additive triples. We show that the set of stable triples defines 
a union of faces of the moment polytope. Moreover these faces may have different dimensions, and many of them have
codimension one. 
\end{abstract}


\section{Introduction}

The complex representation theory of symmetric groups is well understood: the irreducible representations, usually called 
Specht modules, are indexed by partitions and their dimensions are given by the famous hook length formula. But the multiplicative structure of the representation ring has always remained elusive. The multiplicities in tensor products of Specht modules are 
called Kronecker coefficients. They are poorly undestood and notoriously hard to compute. We refer to the introduction to 
\cite{ak1} for a discussion of some of the most basic questions about Kronecker coefficients whose answers remain out of reach. 

That Kronecker coefficients enjoy certain stability properties has been observed by Murnaghan in 1938 \cite{mur1,mur2}.
Such properties are extremely surprising in that they involve representations of differents groups, but they become 
less mysterious once translated in terms of representations of general linear groups, thanks to Schur-Weyl duality.  
 More stability phenomena have been discovered during the last twenty years, and the wealth of examples we are now aware of 
makes more urgent the need to understand and organize them better. That is one of the goals of this note. 

We use the following terminology, taken from \cite{stem} and \cite{ak1}. We denote by $[\lambda]$ the Specht module
associated with the partition $\lambda$. This is an irreducible representation of the symmetric group $S_n$, if $\lambda$ 
is a partition of $n$. Kronecker coefficients are defined by the identity 
$$ [\lambda]\otimes [\mu]=\oplus_\nu g(\lambda,\mu,\nu)[\nu].$$
They are symmetic in $\lambda,\mu,\nu$ and of course, non negative. 

\medskip\noindent {\bf Definition.} 
{\it A triple of partitions $(\lambda,\mu,\nu)$ is}  weakly stable {\it if the Kronecker coefficients} 
$$ g(k\lambda,k\mu,k\nu)=1 \quad  \forall k\ge 1.$$
{\it It is }  stable {\it if $ g(\lambda,\mu,\nu)\ne 0$  and for any triple $(\alpha,\beta,\gamma)$, the sequence of Kronecker 
coefficients $ g(\alpha+k\lambda,\beta+k\mu,\gamma+k\nu)$ is bounded, or equivalently, eventually constant.
We call the asymptotic value of this coefficient a stable Kronecker coefficient.  }

\medskip Stability implies weak stabilty. The converse implication is not known. Conjecturally, the two notions
should be equivalent.

In order to get nice finiteness properties we restrict to partitions whose length, rather that size, are bounded;
the length $\ell(\lambda)$ of a partition $\lambda$ is the number of non zero parts. 
We would then like to understand stability phenomena in relation with the Kronecker semigroup and the 
Kronecker polytohedron. The former is 
$$Kron_{a,b,c}:=\{(\lambda,\mu, \nu),\; \ell(\lambda)\le a,  \ell(\mu)\le b,  \ell(\nu)\le c, \:
g(\lambda,\mu, \nu)\ne 0\}.  $$
This is a finitely generated semigroup. A more precise version of the semigroup property is the elementary, 
but useful {\it  monotonicity property}: 
if  $g(\lambda,\mu, \nu)\ne 0$, then for any triple    $(\alpha,\beta, \gamma)$,
 $$g(\alpha+\lambda,\beta+\mu, \gamma+\nu)\ge g(\alpha,\beta, \gamma).$$

The semigroup $Kron_{a,b,c}$ lives inside a codimension two sublattice of ${\bf Z}^{a+b+c}$, 
because of the obvious condition $|\lambda|=|\mu|=|\nu|$ for a Kronecker coefficient $  g(\lambda,\mu, \nu)$ to be non zero. 
We call this lattice the {\it weight lattice}.
The cone generated by  $Kron_{a,b,c}$  is a rational polyhedral 
cone  $PKron_{a,b,c}$, that we call the Kronecker polyhedron. it is defined by some finite list of linear inequalities, giving the equations of its facets (the maximal faces, of codimension one). 
The number of facets is huge already for small values of the parameters, and certainly grows exponentially with
$a,b,c$ (see \cite{kly} and \cite{vw}). 

In \cite{ak0,ak1}, we showed that certain minimal faces of the Kronecker polyhedron are made of stable triples. 
These minimal faces were defined in terms of certain standard tableaux with the {\it additivity property}.
Let us suppose for simplicity that $c=ab$ (this is not a restriction, since it is well known that for a Kronecker coefficient 
 $  g(\lambda,\mu, \nu)$ 
to be non zero, the condition $\ell(\lambda)\le\ell(\mu)\ell(\nu)$ on the lengths is required). 
Consider a standard   tableau $T$ of rectangular shape $a\times b$. Such a tableau is {\it additive} 
if  there exist increasing sequences 
$x_1<\cdots <x_a$ and $y_1<\cdots <y_b$ of real (or rational) numbers such that 
$$ T(i,j)<T(k,l) \quad\Longleftrightarrow \quad  x_i+y_j<x_k+y_l.$$
The main stability result in \cite{ak1} was the following:

\begin{prop} Let $T$ be any additive standard  tableau of rectangular shape $a\times b$. For any partition 
$\lambda=(\lambda_1, \ldots , \lambda_{ab})$, define two partitions $ a_T(\lambda), b_T(\lambda)$ by 
$$ a_T(\lambda)_i=\sum_{j=1}^b \lambda_{T(i,j)}, \qquad b_T(\lambda)_j=\sum_{i=1}^a \lambda_{T(i,j)}.$$
Then $(\lambda,  a_T(\lambda), b_T(\lambda))$ is a stable triple. 
\end{prop}

Moreover the set of these additive triples, for a fixed $T$, is  exactly the set of lattice points inside a minimal face $f_T$ of $PKron_{a,b,ab}$ defined by this standard tableau. 

We want to stress here that the fact that stable triples can be related to faces of the Kronecker polyhedron is
by no means a surprise. A general statement is the following. 

\begin{prop}
The set $SKron_{a,b,c}$ of weakly stable triples in  $Kron_{a,b,c}$ is the intersection of $Kron_{a,b,c}$ with a union of  faces 
of $PKron_{a,b,c}$.
\end{prop}

More generally, one can associate to any face of the Kronecker polyhedron a positive integer, which gives the order
of growth of the stretched Kronecker coefficients on the interior of the face. This will be discussed in the last section of
this paper. 

Concentrating on (weakly) stable triples, 
it is natural to try to describe the faces of $PKron_{a,b,c}$ that are maximal in $SKron_{a,b,c}$. We will
call the faces of  $PKron_{a,b,c}$ whose intersection with the weight lattice are contained in $SKron_{a,b,c}$ 
the  {\it stable faces}, and those that are 
maximal in  $SKron_{a,b,c}$, the {\it maximal stable faces} (which we don't expect a priori to be maximal faces in 
$PKron_{a,b,c}$, or facets). Among 
many other questions, we can ask: what is the maximal dimension of such a face? what can be their dimensions? 
could they all be of the same dimension? can the additive stable faces be maximal in $SKron_{a,b,c}$? more
generally, what are the stable faces containing a given additive stable face? 

The main goal of this note is to  answer some of these questions, in particular the last one, 
and draw some unexpected consequences. 
In \cite{ak1} we explained how to describe the local structure of the Kronecker polyhedron around an additive face. 
Among the faces that contain such an additive face, we will distinguish those that have a property that we will call
{\it strong simpliciality}. We will prove:

\begin{theo}
Among the faces of  $PKron_{a,b,c}$ that contain an additive face, the stable ones are exactly those 
that are strongly simplicial.  
\end{theo}

A priori, we would have expected the stable faces to be very special, in particular to have high codimension.
Surprisingly, our Theorem has the following consequence:

\begin{coro}
The polyhedral cone  $PKron_{a,b,c}$ always contains   stable facets.
\end{coro}

This means that there exist families of stable triples of the largest possible dimension. 
It would be extremely interesting to 
have a full classification. We can give
many explicit examples of strongly simplicial facets and show that there always exist many of them (Proposition 4).
 We can also describe their structure, which is that of a cone over hypercube (Proposition 3). The vertices of
this hypercube are in bijection with the additive faces contained in the facet.  

Another striking phenomenon is the following. Consider an additive face, and the maximal stable faces that 
contain it. It may very well happen that these maximal faces have different dimensions! In fact it seems quite
plausible that the maximal stable faces can have all the possible dimensions between the smallest and maximal
possible dimensions. In particular the set of stable triples seems to have a very intricate structure in general. 

\section{Strongly simplicial faces}

\subsection{The geometric method}

Let us briefly recall the main features of the geometric method used in \cite{ak0,ak1} in order to approach Kronecker coefficients. 
Let $A, B$ be complex vector spaces of finite dimensions $a, b$.  By Schur-Weyl duality, Kronecker coefficients are the 
multiplicities of the Schur powers $S_\lambda (A\otimes B)$, when decomposed into irreducible representations for $GL(A)
\times GL(B)$. By the Borel-Weil theorem, 
$$ S_\lambda (A\otimes B)=H^0(Fl (A\otimes B), L_\lambda)$$ 
for a suitable linearized line bundle  $L_\lambda$ on the variety $Fl (A\otimes B)$ of complete flags in $A\otimes B$.  
A standard tableau $T$ defines an embedding 
$$ \iota_T : Fl(A)\times Fl(B)\hookrightarrow Fl(A\otimes B).$$
The induced map on equivariant Picard groups 
$$\iota_T^* : Pic( Fl(A\otimes B))\simeq {\bf Z}^{ab}\rightarrow Pic( Fl(A)\times Fl(B))\simeq {\bf Z}^a\times {\bf Z}^b$$ 
is precisely our map $\lambda\mapsto (a_T(\lambda),
b_T(\lambda))$ when expressed in natural basis. In particular, restriction gives a nonzero map 
\begin{equation}\label{res}
H^0(Fl (A\otimes B), L_\lambda)\longrightarrow H^0(Fl (A)\times Fl(B), L_{a_T(\lambda)}\otimes L_{b_T(\lambda)})
=S_{a_T(\lambda)}A\otimes S_{b_T(\lambda)}B,
\end{equation}
implying that the Kronecker coefficient $g(\lambda, a_T(\lambda), b_T(\lambda))$ is positive. Then we can  define a filtration of 
$H^0(Fl (A\otimes B), L_\lambda)$ by the order of vanishing on $Fl (A)\times Fl(B)$. This allows to define an injective map
\begin{equation}\label{taylor}
H^0(Fl (A\otimes B), L_\lambda)\hookrightarrow H^0(Fl (A)\times Fl(B), L_{a_T(\lambda)}\otimes L_{b_T(\lambda)}
\otimes S^*N^*),
\end{equation}
where $N$ denotes the normal bundle of the embedding $\iota_T$, and  $ S^*N^*$ is the symmetric algebra of the dual 
bundle, the conormal bundle. This map must be thought of as taking a section of 
$L_\lambda$, to its Taylor expansion in the normal directions to  $Fl (A)\times Fl(B)$.
Moreover, if $\lambda$ is strictly decreasing, the line bundle 
$L_\lambda$ is very ample. By the usual properties of ample bundles, the previous map becomes surjective onto every finite
part of $ S^*N^*$  if $L_\lambda$ is sufficiently ample (that is, if the differences $\lambda_i-\lambda_{i+1}$ are large enough). 
This shows that the multiplicities in $ S_\lambda (A\otimes B)=H^0(Fl (A\otimes B), L_\lambda)$,
otherwise said,  the Kronecker coefficients, are somehow controlled by the normal bundle. 

This works particularly well when the embedding $\iota_T$ is {\it convex}, in the sense that the weights of the normal bundle are
contained in a strictly convex cone. Combinatorially, this exactly means that the tableau $T$ is additive. Then the Kronecker
coefficient $g(\alpha+k\lambda, \beta+ka_T(\lambda), \gamma+kb_T(\lambda))$ is bounded by the multiplicity of 
 $L_{\beta-a_T(\alpha)}\otimes L_{\gamma-b_T(\alpha)}$ inside $ S^*N^*$, and the latter is finite by convexity. This implies
that we can focus on a finite part of this algebra, independantly of $k$. But then, if $\lambda$ is strict and $k$ is large 
enough, the surjectivity of (\ref{taylor}) in finite degrees implies that we have equality between the latter multiplicity, and the Kronecker coefficient. In particular this coefficient does not depend on $k$, when big enough: this is the stability 
phenomenon. But of course we get
much more information since we are in principle able to compute the stable Kronecker coefficients, directly from the normal bundle. 

Combinatorially, the weights of the conormal bundle are determined as follows. Denote by $e_1,\ldots ,e_a$ and $f_1,\ldots ,f_b$
basis of the character lattices of maximal tori in $GL(A)$ and $GL(B)$. If $T(i,j)=k$ (ie the box $(i,j)$ is numbered $k$ in $T$), 
let $g_k=e_i+f_j$. Then the weights of the normal bundle are the differences $g_\ell-g_k$ for $\ell>k$. Among those 
weights, the horizontal and vertical ones are those of the form $e_p-e_q$ and $f_p-f_q$. They will appear repeatedly, 
in fact in the conormal bundle their multiplicities are  $b-1$ for the horizontal ones, $a-1$ for the vertical ones. In 
particular these multiplicities will be bigger than one as soon as we suppose that $a,b>2$. All the other weights have 
multiplicity one. Of course, the multiplicity of $L_{\beta-a_T(\alpha)}\otimes L_{\gamma-b_T(\alpha)}$ inside $ S^*N^*$,
which gives the stable Kronecker coefficient, can be obtained as the number of ways to express the weight 
$(\beta-a_T(\alpha),\gamma-b_T(\alpha))$ as a non negative linear combination of the weights of the conormal bundle, 
considered with their multiplicities. 

Of course this is possible only when $(\beta-a_T(\alpha),\gamma-b_T(\alpha))$ belongs to the cone generated by those
weights, which we call the {\it conormal cone}.  This cone gives a local picture of the Kronecker polyhedron locally
around $f_T$. In particular any face of the latter containing $f_T$, can be identified with a face of the conormal
cone, and conversely. 

\subsection{Perturbations of additive triples}

Now consider a triple of the form $(\lambda, a_T(\lambda)+\sigma, b_T(\lambda)+\theta),$ where $\sigma$ and $\theta$
are not necessarily partitions, but sequences (or weights) such that $ a_T(\lambda)+\sigma$ and $b_T(\lambda)+\theta$ are partitions.

By the injectivity of (\ref{taylor}), the Kronecker coefficient $g(k\lambda, k(a_T(\lambda)+\sigma), k(b_T(\lambda)+\theta))$
is bounded by the multiplicity of $(k\sigma, k\theta)$ as a weight of  $ S^*N^*$. If we suppose that  the line generated by
$(\sigma, \theta)$ belongs to the conormal cone, this multiplicity will eventually become positive, and we expect it
to grow to infinity with $k$. But this is not necessarily the case: the multiplicity will remain bounded if $(\sigma, \theta)$
belongs to a face of the cone which is  {\it strongly simplicial}. By this we mean 
that the weights of the conormal bundle contained in the face, 
considered with their multiplicities, define a basis of the linear space  generated by the face. Then the multiplicity will be 
$0$ or $1$, the second possibility occuring exactly when $(k\sigma, k\theta)$ belongs to the lattice generated by the latter weights. 

Let us insist on the definition of strongly simplicial faces. 

\medskip\noindent 
{\bf Definition}. {\it A face $F$ of the Kronecker polytope is} strongly simplicial {\it  if it contains an additive face $f_T$ such 
that locally around $f_T$, the face $F$ corresponds to a face in the cone generated by the conormal bundle which
is strongly simplicial in the sense that:
\begin{enumerate}
\item it is a face of dimension $d$ generated by $d$ vectors  $g_{k_1+1}-g_{k_1}, \ldots , g_{k_d+1}-g_{k_d}$, 
\item none of these vectors is horizontal (unless $b=2$) or vertical (unless $a=2$),
\item no other vector of the form $g_p-g_q$ belongs to the face, 
\item in particular the pairs $(k_1,k_1+1), \ldots , (k_d,k_d+1)$ do not intersect.
\end{enumerate} }

The structure of a strongly simplicial face is not difficult to describe. Recall that an additive tableau $T$ is defined by parameters 
$x_1<\cdots <x_a$ and $y_1<\cdots <y_b$ such that $ T(i,j)<T(k,l)$ if and only if  $x_i+y_j<x_k+y_l.$ Of course
these parameters are not unique. In fact the tableau $T$ really corresponds to a connected component $C_T$ of the complement
of the collection of hyperplanes defined by the equations $x_i+y_j=x_k+y_l$ inside the parameter space.

 Locally around the additive face $f_T$, the Kronecker
polyhedron is, by hypothesis, the simplicial cone over the vectors  $g_{k_1+1}-g_{k_1}, \ldots , g_{k_d+1}-g_{k_d}$. 
Let us choose one of them,  say $g_{k_s+1}-g_{k_s}$. Since it is neither horizontal nor vertical, we can exchange the 
entries $k_s$ and $k_s+1$ in $T$ and obtain another standard tableau $T_s$. We claim that $T_s$ is again additive. 
Indeed, if the entries  $k_s$ and $k_s+1$ of $T$ appear in boxes $(i,j)$ and $(k,l)$, the fact that 
 $g_{k_s+1}-g_{k_s}$ is an extremal vector of the cone implies that the hyperplane  $x_i+y_j=x_k+y_l$
is really a facet of $C_T$. Crossing this facet we get into a  component corresponding to $T_s$, which is therefore additive.

Iterating the process, we deduce that the $2^d$ standard tableaux obtained by considering all the possibilities to
exchange the entries   $(k_1,k_1+1) \ldots (k_d,k_d+1)$, are all additive. Moreover the Kronecker polyhedron,
around each of the corresponding additive faces, is described by the same cone, up to a change of signs for the
generators. This implies that our strongly simplicial face is contained in the set of triples 
\begin{equation}\label{sf}
(\lambda,\mu,\nu)=(\lambda, a_T(\lambda), b_T(\lambda))+\sum_{i=1}^du_s(g_{k_s+1}-g_{k_s}),
\end{equation}
with $0\le u_s\le \lambda_{k_s}-\lambda_{k_s+1}$. Note that the coefficients $u_1,\ldots ,u_d$ need to be integers. 
But it is a priori possible that we get a triple of partitions $(\lambda,\mu,\nu)$ given by the same expression but with 
rational coefficients  $u_1,\ldots ,u_d$, not all integral. In this case, the Kronecker coefficient $g(\lambda,\mu,\nu)$
would certainly be zero. 

Otherwise said, the identity (\ref{sf}) defines a lattice $L_F$, which could be a proper sublattice of the intersection of the 
weight lattice with the linear span of $F$.  In this lattice, $F$ is simply defined by the inequalities $0\le u_s\le \lambda_{k_s}-\lambda_{k_s+1}$ for $1\le s\le d$. Recall that $d$ is the number of generators of the face in the normal directions
of an additive face it contains. In  particular the codimension of $F$ is the codimension of an additive face (that is, $a+b-2$)
minus $d$. We get the following description of strongly simplicial faces. 

\begin{prop}
A strongly simplicial face $F$ of codimension $\delta$ in the Kronecker polyhedron is a cone over a hypercube of dimension $a+b-2-\delta$.   
\end{prop}

The main result of this paper is the following. 

\begin{theo}
A  strongly simplicial face $F$ of the Kronecker polyhedron is stable. More precisely, any point in $F$ is a stable triple 
if it belongs to $L_F$, and the corresponding Kronecker coefficient is zero otherwise. 
\end{theo}

\noindent {\it Proof.}
Consider a triple of the form $(k\lambda+\alpha, k(a_T(\lambda)+\sigma)+\beta, k(b_T(\lambda)+\theta))+\gamma),$ 
where as before $(\sigma, \theta)$
belongs to the simplicial face corresponding to $F$ in the conormal cone of the additive face $f_T$. 
As we have seen, the corresponding Kronecker coefficient is bounded by the multiplicity of 
$(k\sigma+\beta-a_T(\alpha), k\theta+\gamma-b_T(\alpha))$ as a weight of  $ S^*N^*$.
Suppose that we have expressed this weight as a non negative integer linear combination 
$t_1\eta_1+\cdots +t_N\eta_N$ of the weights $\eta_1, \ldots ,\eta_N$ of the conormal bundle,
considered with their multiplicities.  Suppose these weights are indexed in such a way that the first $d$ 
generate our simplicial face. By projection along the direction of this face, we get a relation of the form
\begin{equation}\label{proj}
p_F(\beta-a_T(\alpha),\gamma-b_T(\alpha))=t_{d+1}p_F(\eta_{d+1})+\cdots + t_{N}p_F(\eta_{N}),
\end{equation}
where $p_F$ denotes the projection. But the projected weights  $p_F(\eta_{d+1}), \ldots , p_F(\eta_{N})$
generate a strictly convex cone, so the latter equation has only finitely many non negative integer solutions 
$(t_{d+1}, \ldots , t_N)$. These solutions do not depend on $k$, and for each of these, the original equation 
has at most one solution in $(t_1,\ldots , t_d)$, since it can be considered as an equation in the simplicial face 
$F$. This proves that the Kronecker coefficient  $g(k\lambda+\alpha, k(a_T(\lambda)+\sigma)+\beta, k(b_T(\lambda)+\theta))+\gamma),$ is bounded independently of $k$. 

This is precisely the definition of stability, up to the fact that the Kronecker coefficient $g(\lambda, a_T(\lambda)+\sigma, b_T(\lambda)+\theta)$ must be equal to one. Recall that by \cite{stem}, Propoition 3.2, the only alternative is 
that it is equal to zero. So what remains to prove is that if $(\lambda, \mu, \nu)$ is a point of $F$ that also belongs 
to the lattice $L_F$, the Kronecker coefficient  $g(\lambda, \mu, \nu)$ cannot be zero. 

To check this we will use that  $(\lambda, \mu, \nu)$ is given by (\ref{sf}) for some integer coefficients 
 $u_1,\ldots ,u_d$ such that $0\le u_s\le \lambda_{k_s}-\lambda_{k_s+1}$ for all $s$. Denote by 
$\omega_t$ the partition of $t$ with $t$ parts equal to one. Recall that we denoted by $T_s$ the standard 
tableau obtained by exchanging the entries $k_s$ and $k_s+1$ in $T$. It is straightforward to check that 
$$(\omega_{k_s}, a_{T_s}(\omega_{k_s}), b_{T_s}(\omega_{k_s}))-
 (\omega_{k_s}, a_T(\omega_{k_s}), b_T(\omega_{k_s}))= g_{k_s+1}-g_{k_s}.$$
This allows to rewrite   (\ref{sf}) as 
$$(\lambda,\mu,\nu)=
(\theta, a_T(\theta), b_T(\theta))+\sum_{s=1}^du_s(\omega_{k_s}, a_{T_s}(\omega_{k_s}), b_{T_s}(\omega_{k_s})),$$
where $\theta=\lambda -\sum_{s=1}^du_s\omega_{k_s}$. Since $u_s\le \lambda_{k_s}-\lambda_{k_s+1}$ for all $s$,
this $\theta$ is again a partition. 
Since $T$ and the $T_s$ are all additive, we know that $g(\theta, a_T(\theta), b_T(\theta))=1$ and 
$g(\omega_{k_s}, a_{T_s}(\omega_{k_s}), b_{T_s}(\omega_{k_s}))=1$ for all $s$. In particular all these Kronecker coefficients are non zero, and from the semigroup property we deduce that  $g(\lambda,\mu,\nu)$ is positive. $\Box$

\begin{coro}
A strongly simplicial face is the non negative integral span of the additive faces it contains. 
Moreover any additive face  is properly contained in some strongly simplicial face.  
\end{coro} 

\noindent {\it Proof.} The first statement means that any stable triple in $F$ can be obtained as a linear combination with positive integer coefficients, of some stable triples in the additive faces contained in $F$. This is what we established in the proof of the Theorem. 

 For the second statement, simply observe that at least one face of the Kronecker polyhedron, that contains $f_T$ and has 
dimension one more, must be simplicial. Indeed, these faces correspond to the minimal generators of the
conormal cone, and they are simplicial exactly for those generators that are neither horizontal nor vertical.
But the generators cannot be all horizontal or vertical, since otherwise inside $T$, the 
integer $k+1$ would always be  South-East of $k$, which is absurd.   $\Box$

\medskip\noindent  {\it Remarks.}  One can wonder if there can be non trivial arithmetic conditions on the strongly simplicial faces,
for the Kronecker coefficients to be non zero? This would mean that $L_F$ is really a proper sublattice of the intersection of
$F$ with the weight lattice.  
This seems a priori possible but we have no example of such a phenomenon. 

One can also wonder if  stable Kronecker coeficients, when one considers strongly simplicial faces,  count points in some polytopes,
as they do on additive faces \cite{ak1}. In the proof above we indeed bounded the stretched Kronecker coefficients by 
numbers of points in some polytopes, but it is not clear that this bound coincides with the stable Kronecker coefficient. In the 
additive case this follows from an ampleness argument, which does not apply in this more general situation.
  
\subsection{Strongly simplicial facets}

In \cite{ak1} we gave a combinatorial description of the facets of the Kronecker polytope containing a given 
additive face $f_T$. These facets are in bijection with the {\it maximal relaxations} compatible with $T$, where 
a maximal relaxation $R$ is given by an additive (non standard) tableau defined by sequences
 $x_1\le \cdots \le x_a$  and  $y_1\le \cdots \le y_b$ such that the sums $R(i,j)=x_i+y_j$ are not
necessarily distinct. What we require is that the set of vectors $e_i+f_j-e_k-f_l$, for    $R(i,j)=R(k,l)$,
has maximal rank $r=a+b-3$. Such a family of vectors being given,  the  sequences $x_1\le \cdots \le x_a$  
and  $y_1\le \cdots \le y_b$ are uniquely defined up to translation, and multiplication by the same positive 
number.  It is convenient to define the tableau $R$ uniquely by letting $x_1=y_1=0$, and asking the 
two sequences to be made of integers, with no common divisor.  The compatibility condition with a standard
tableau $T$ is that  $R(i,j)<R(k,l)$ implies  $T(i,j)<T(k,l)$. Otherwise said, $R$ defines a partial order on the 
boxes in the rectangle $a\times b$, which is refined by the total order defined by $T$. The equation of the 
facet $F_R$ is then  given by 
$$\sum_{i=1}^ax_i\mu_i+\sum_{j=1}^by_j\nu_j=\sum_{i=1}^a\sum_{j=1}^b(x_i+y_j)\lambda_{T(i,j)},$$
where $T$ is any standard tableau compatible with $R$.
 
Can such a maximal relaxation $R$ define a strongly  simplicial facet? This would mean that $R$ is defined by strictly 
increasing sequences, and that there exists exactly $r=a+b-3$ values of $R$ appearing twice, the 
corresponding difference vectors being independent. In terms of the hyperplanes of equations 
$x_i+y_j-x_k-y_l=0$, and the arrangement they define in the open cone defined by  
$0=x_1< \cdots < x_a$  and  $0=y_1< \cdots < y_b$, such an $R$ corresponds to a point 
where exactly $r$ hyperplanes meet transversaly. Recall that this transversality property implies that 
any of the $2^r$ standard tableaux $T$ compatible with $R$ is additive.

\medskip Another unexpected fact is that in general, 
there exist surprisingly many strongly simplicial facets!

\begin{prop}
$PKron(a,b,ab)$ contains at least ${a+b-4 \choose b-2}$ strongly simplicial facets.
\end{prop} 

\noindent {\it Proof.} One can construct tableaux defining strongly simplicial facets by a simple induction: 
suppose that a tableau $S$ defines a simplicial facet
for the format $a\times (b-1)$. Then we get one for the format $a\times b$ by adding a column  defined by 
$y_b=x_a+y_{b-1}$. Of course this also works for rows. So starting from the tableau defining the 
unique simplicial face in format $2\times 2$, we can construct    ${a+b-4 \choose b-2}$ 
strongly simplicial facets in format $a\times b$ by choosing to 
apply the previous process on rows or columns successively,  in all possible orders. $\Box$  

\subsection{Examples}

Let us examine in more detail the low dimension cases. 

\medskip\noindent {\bf Example 1.} For $a=b=2$ there is exactly one additive face (up to symmetry). This additive 
face is the intersection of two facets, one of which is strongly simplicial. On the additive face we get  
$$g( (\lambda_1,\lambda_2, \lambda_3, \lambda_ 4), (\lambda_1+\lambda_2,  \lambda_3+\lambda_ 4),
(\lambda_1+\lambda_3, \lambda_2+ \lambda_ 4) )=1,$$
and for the strongly simplicial facet we get the more general statement that 
$$g( (\lambda_1,\lambda_2, \lambda_3, \lambda_ 4), (\mu_1,  \mu_2), (\nu_1,  \nu_2))=1$$
when $\mu_1-\nu_2=\lambda_1-\lambda_4$ and $\lambda_1+\lambda_3\le\mu_1\le\lambda_1+\lambda_2$.
Moreover all these triples are stable.

\medskip\noindent {\bf Example 2.} For $a=b=3$ there exist $42$ standard 
tableaux fitting in a square of size three, 
among which $36$ are additive. The number of maximal relaxations  is $17$. They are encoded in the following tableaux:
$$\begin{array}{ccccc}
F_1^+ \begin{array}{ccc}
0&0&0 \\ 0&0&0 \\ 1&1&1
\end{array}\;\; &
F_2^+\begin{array}{ccc}
0&0&0 \\ 1&1&1 \\ 1&1&1
\end{array}\;\; &
F_3^+\begin{array}{ccc}
0&0&1 \\ 1&1&2 \\ 2&2&3
\end{array}\;\; &
F_4^+\begin{array}{ccc}
0&1&1 \\ 1&2&2 \\ 2&3&3
\end{array}\;\; &
F_5^+\begin{array}{ccc}
0&1&2 \\ 2&3&4 \\ 3&4&5
\end{array} \\
 &&&& \\
F_6\begin{array}{ccc}
0&0&1 \\ 0&0&1 \\ 1&1&2
\end{array}\;\; &
F_7^+\begin{array}{ccc}
0&0&1\\ 1&1&2 \\ 1&1&2
\end{array}\;\; &
F_8\begin{array}{ccc}
0&1&1 \\ 1&2&2 \\ 1&2&2
\end{array}\;\; &
F_9\begin{array}{ccc}
0&1&2 \\ 1&2&3 \\ 2&3&4
\end{array}\;\; &
F_{10}^+\begin{array}{ccc}
0&1&2 \\ 1&2&3 \\ 3&4&5
\end{array}
\end{array}$$
\noindent and for each tableau $F_i^+$ there is another one denoted $F_i^-$ obtained by diagonal symmetry. 

Recall that additive faces have dimension four. A detailed analysis yields the following result: 

\begin{prop}
For $a=b=3$, the maximal strongly simplicial faces are, up to diagonal symmetry:
\begin{enumerate}
\item in codimension one, $F_5^+$ and $F_{10}^+$,
\item in codimension two,  $F_3^+\cap F_{4}^+$, $F_3^+\cap F_9$, $F_4^+\cap F_9$,
\item in codimension three, $F_6\cap F_7^+\cap F_9$ and  $F_7^+\cap F_8\cap F_9$.
\end{enumerate}
\end{prop}

Let us describe the sets of triples $(\lambda,\mu,\nu)$ on these strongly simplicial faces. We will use the notations
$\lambda_{ij}=\lambda_{i}+\lambda_{j}$ and $\lambda_{ijk}=\lambda_{i}+\lambda_{j}+\lambda_{k}$. 

\medskip\noindent {\bf $F_5^+$} is defined by the equation $2\mu_2+3\mu_3+\nu_2+2\nu_3=
\lambda_2+2\lambda_3+2\lambda_4+3\lambda_5+3\lambda_6+4\lambda_7+4\lambda_8+5\lambda_9$
and the   inequalities 
$ \lambda_{124}\le\mu_1\le \lambda_{123},$
$ \lambda_{123}+ \lambda_{146}\le\mu_1+\nu_1\le \lambda_{123}+\lambda_{145},$
$\lambda_{12}-\lambda_{79}\le \mu_1-\nu_3\le  \lambda_{12}-\lambda_{89}.$

\medskip\noindent {\bf $F_{10}^+$}  is defined by the equation $\mu_2+3\mu_3+\nu_2+2\nu_3=
\lambda_2+\lambda_3+2\lambda_4+2\lambda_5+3\lambda_6+3\lambda_7+4\lambda_8+5\lambda_9$
and the inequalities
$\lambda_{13}-\lambda_{89}\le \nu_1-\mu_3\le  \lambda_{12}-\lambda_{89},$
$ \lambda_{569}+ \lambda_{789}\le\mu_3+\nu_3\le \lambda_{469}+\lambda_{789},$
$ \lambda_{789}\le\mu_3\le \lambda_{689}.$

\medskip\noindent {\bf  $F_3^+\cap F_{4}^+$} is defined by the two equations $\nu_2=\lambda_{258}$ and 
$\mu_2+2\mu_3+\nu_3=\lambda_{369}+\lambda_{456}+2\lambda_{789}$ and the inequalities
$\lambda_{124}\le \mu_1\le \lambda_{123}$, $ \lambda_{789}\le \mu_3\le \lambda_{689}.$

\medskip\noindent {\bf  $F_3^+\cap F_{5}$} is defined by the two equations $\nu_1=\lambda_{136}$ and 
$\mu_1-\mu_3+\nu_2=\lambda_{124}+\lambda_{258}-\lambda_{689}$ and the  inequalities
$\lambda_{125}\le \mu_1\le \lambda_{124}$, $\lambda_{689}\le \mu_3\le \lambda_{679}.$

\medskip\noindent {\bf  $F_4^+\cap F_{5}$} is defined by the two equations $\nu_3=\lambda_{479}$ and 
$\mu_1-\mu_3-\nu_2=\lambda_{124}-\lambda_{258}-\lambda_{689}$ and the inequalities
$\lambda_{134}\le \mu_1\le \lambda_{124}$, $\lambda_{689}\le \mu_3\le \lambda_{589}.$

\medskip\noindent {\bf  $F_6\cap F_7^+\cap F_9$}  is defined by the three equalities  $\mu_1+\nu_1=\lambda_{125}
+\lambda_{136}$,  $\mu_2=\lambda_{348}$,  $\nu_2=\lambda_{247}$ and the inequalities  $\lambda_{126}\le \mu_1\le \lambda_{125}.$

\medskip\noindent {\bf  $F_7^+\cap F_8\cap F_9$} is defined by the three equalities  $\mu_1+\nu_1=\lambda_{124}
+\lambda_{135}$,  $\mu_2=\lambda_{356}$,  $\nu_2=\lambda_{267}$ and the inequalities  $\lambda_{125}\le \mu_1\le \lambda_{124}.$

\medskip
There are no arithmetic constraints on these stongly simplicial faces, 
so the Kronecker coefficients are always equal to one and all these
triples are stable. 

\smallskip
Note also that the additive face defined by the standard tableau $$ \begin{array}{ccc}
1&2&3 \\ 4&5&7 \\ 6&8&9
\end{array}$$
is contained in both $F_5^+$ and   $F_3^+\cap F_{4}^+$, showing that 
an additive face can be contained in two maximal strongly simplicial faces of different dimensions! 
This indicates that the structure of the set of additive triples must be quite intricate in general. 

\medskip\noindent {\bf Example 3.} For $a=b=4$ there are $6660$ additive tableaux and $457$ maximal relaxations, according to 
\cite{kly,vw}. Among these, we know $43$ strongly simplicial ones, among which:

$$\begin{array}{cccc}

0&1&2&3 \\  1&2&3&4 \\ 4&5&6&7 \\7&8&9&A
\end{array} 
\qquad 
\begin{array}{cccc}
0&1&2&5 \\ 1&2&3&6 \\ 3&4&5&8 \\ 7&8&9&C
\end{array} \qquad 
\begin{array}{cccc}
0&1&2&5 \\ 1&2&3&6 \\ 3&4&5&8 \\ 8&9&A&D
\end{array} \qquad 
\begin{array}{cccc}
0&1&2&5 \\ 2&3&4&7 \\ 3&4&5&8 \\ 8&9&A&D
\end{array} \qquad
\begin{array}{cccc}
0&1&2&6 \\ 2&3&4&8 \\ 3&4&5&9 \\ 5&6&7&B
\end{array}$$

 $$
\begin{array}{cccc}
0&1&2&7 \\ 1&2&3&8 \\ 3&4&5&A \\ 5&6&7&C
\end{array} \qquad 
\begin{array}{cccc}
0&1&2&7 \\ 2&3&4&9 \\ 3&4&5&A \\ 5&6&7&C
\end{array} \qquad 
\begin{array}{cccc}
0&1&4&6 \\ 2&3&6&8 \\ 3&4&7&9 \\ 7&8&B&D
\end{array} \qquad
\begin{array}{cccc}
0&1&5&7 \\ 5&6&A&C \\ 6&7&B&D \\ B&C&G&I
\end{array} \qquad 
\begin{array}{cccc}
0&2&3&6 \\ 3&5&6&9 \\ 5&7&8&B \\ 7&9&A&D
\end{array}$$

(The symbols $A,B$ and so on stand for $10, 11$ and so on.)
It would be interesting to have the complete list.

\subsection{Symmetries}

There exist two natural involutions on the set of additive tableaux. Recall that an additive tableau can be 
defined by increasing sequences 
$x_1<\cdots <x_a$ and $y_1<\cdots <y_b$ such that the sums $x_i+y_j$ are distinct. Then we can replace each of this 
sequence by the opposite one, reordered increasingly. Since this preserves the set of hyperplanes of equations
$x_i+y_j=x_k+y_l$, this defines two commuting involutions on the set of additive tableaux, and then also on the 
set of maximal relaxations, and on the subset of simplicial relaxations. 

\section{Non simplicial faces}

\subsection{The degree of a face}

Recall that a stretched Kronecker coefficient $g(k\lambda,k\mu,k\nu)$ is quasipolynomial: there exists  a collection 
of polynomials $P_0, \ldots ,P_{p-1}$, such that 
$$g(k\lambda,k\mu,k\nu)=P_i(k) \qquad \mathrm{for}\quad k=i\; (\mathrm{mod}\; p).$$ 
By the monotonicity property,  $P_{i+j}(k+\ell)\ge P_i(k)$ as soon as $P_j(\ell)\ne 0$. This implies that among the 
polynomials $P_0, \ldots ,P_{p-1}$, those that are not identically zero have the same degree $d$, and the same leading term
as well. We call $d=d(\lambda,\mu,\nu)$ the  degree of the triple $(\lambda,\mu,\nu)$. 
For example weakly stable triples have degree zero, and a triple of degree zero is one that has a weakly stable multiple. 

Another straightforward consequence of the monotonicity property, and of the convexity of the faces, is the 
following statement. 

\begin{prop}\label{degree}
Let $F$ be a face of the Kronecker polyhedron. The degree is constant on the interior of $F$, and can only decrease, or 
remain the same, on its boundary faces. 
\end{prop}

\noindent {\bf Definition.}  {\it Let $F$ be a face of the Kronecker polyhedron. We define its}  degree {\it as the degree of its interior 
points.  For example, any additive face, more generally any strongly simplicial face, has degree zero.}

\subsection{The defect of simpliciality}

For a non simplicial face, containing an additive face, we will show that the degree can be read off directly 
on the normal bundle. 

\medskip\noindent {\bf Definition.}  {\it A face $F$ of the Kronecker
polyhedron will be called}  $\delta$-simplicial {\it if there exists an additive face $f_T$ in $F$, such that the face $f$ of the conormal
cone corresponding to $F$ is $\delta$-simplicial. By this we mean that $f$ has dimension $d$, but contains $d+\delta$ weights 
of the conormal bundle, counted with their multiplicity.}

\medskip Strongly simplicial is therefore the same as $0$-simplicial. Note also that starting from a face   $F$ of the Kronecker
polyhedron, the integer  $\delta$ will not depend on the minimal face $f_T$ contained in $F$. This is a consequence of the 
following statement:

\begin{theo}
A $\delta$-simplicial face $F$ of the Kronecker polyhedron has degree $\delta$.
\end{theo}

\noindent {\it Proof.} We consider $F$ with the additive face $f_T$, and we identify $F$ with the corresponding face of 
the conormal cone. We consider stretched Kronecker coefficients $g(k\lambda, k(a_T(\lambda)+\sigma), k(b_T(\lambda)+\theta))$, 
where the weight $(\sigma, \theta)$ belongs to the linear span of the face. Denote this Kronecker coefficient by $g_k$. It may a priori happen that the lattice generated by the weights of the conormal bundle belonging to the face does not contain $(\sigma, \theta)$. In general there exists a minimal
integer $e$, depending on $(\sigma, \theta)$, such that $e(\sigma, \theta)$ belongs to this lattice. If $k$ is not divisible by $e$, 
then $k(\sigma, \theta)$ does not belong to the lattice and $g_k=0$. If we restrict to those $k$ that are divisible by $e$, 
then the number of ways to express $k(\sigma, \theta)$ as a non negative integer linear combination of weights of the 
conormal bundle certainly grows like $k^\delta$. By the injectivity of (\ref{taylor}), this implies that the growth of $g_k$ is 
at most in   $k^\delta$.

To get to the required conclusion, we must control the surjectivity of (\ref{taylor}). The key point is the following 
general statement. 

\begin{lemm}
Let $L$ be an ample line bundle, $M$ a globally generated line bundle on a smooth complex projective variety $X$. 
Let $\iota: Y\hookrightarrow X$ be the embedding of a smooth subvariety, and denote by $N$ the normal bundle. 
Then  there exists integers $m$ and $n$, not depending on $M$, such that the natural map 
$$H^0(X,I_Y^d\otimes L^a\otimes M)\longrightarrow H^0(Y, S^dN^*\otimes \iota^*( L^a\otimes M))$$
is surjective when $a\ge md+n$.   
\end{lemm} 

If we apply this statement to $\iota_T$, we deduce that there exists  integers $m_T, n_T$ such that (\ref{taylor}) 
is surjective up to degree $d$ as soon as $\lambda_i-\lambda_{i+1}\ge m_Td+n_T$ for each $i$. 
Replacing $\lambda$ by $k\lambda$ we get the surjectivity up to degree  $(k-n_T)/m_T$. This yields a lower bound 
for $g_k$ of order $(k-n_T/m_T)^\delta$, and the claim follows. $\Box$

\medskip \noindent {\it Proof of the Lemma.} To get the surjectivity it is enough to prove that 
$$  H^1(X,I_Y^{d+1}\otimes L^a\otimes M)=0. $$
 Let $\pi : Z\rightarrow X$ be the blow-up of $Y$, and $E$ the exceptional divisor. Since $\pi_*O_Z(-iE)=I_Y^i$
and there are no higher direct images, we are reduced to proving that 
$$  H^1(Z,O_Z(-(d+1)E)\otimes \pi^*(L^a\otimes M))=0. $$
The canonical line bundle of $Z$ is $K_Z=\pi^*K_X\otimes O_Z((c-1)E)$, if $c$ denotes the codimension of $Y$ in $X$. 
So we can rewrite the previous condition as 
$$  H^1(Z,K_Z\otimes O_Z(-(d+c)E)\otimes \pi^*(L^a\otimes M\otimes K_X^{-1}))=0. $$
 We can find an $a_0$ such that $L^{a_0}\otimes K_X^{-1}$ is ample. Moreover the exists $b_0$ such that  $I_Y\otimes L^{b_0}$ is generated by global sections, hence also $O_Z(-E)\otimes \pi^*L^{b_0}$.  Then    $O_Z(-(d+c)E)\otimes \pi^*(L^a\otimes M\otimes K_X^{-1})$ is nef and big as soon as  
$a\ge a_0+b_0(d+c)$, and the required vanishing follows from the Kawamata-Viehweg vanishing theorem. $\Box$  
 
\medskip
\noindent {\it Acknowledgements.} This paper was  in Montr\'eal at the Centre de Recherches Math\'ematiques (Universit\'e de Montr\'eal) and the CIRGET  (UQAM). The author warmly thanks these institutions for their generous hospitality. We also thank 
Mateusz Michalek for his help regarding the combinatorics of simplicial facets. In particular Proposition 6 is due to him.

\end{document}